\documentclass[11pt]{article}
\usepackage[T2A]{fontenc}
\usepackage{amsmath}
\usepackage{amssymb}
\usepackage{arydshln}
\def\a{\mathbf{a}}
\def\hh{\mathbf{h}}
\def\pp{\mathbf{p}}

\def\m{\mathbf{m}}
\def\f{\mathbf{f}}

\def\b{\mathbf{b}}
\def\N{\mathbb{N}}
\def\Z{\mathbb{Z}}

\def\C{\mathbb{C}}

\def\llambda{\boldsymbol{\lambda}}
\def\ggamma{\boldsymbol{\gamma}}
\def\aalpha{\boldsymbol{\alpha}}
\def\bbeta{\boldsymbol{\beta}}
\def\llambda{\boldsymbol{\lambda}}

\def\ddelta{\boldsymbol{\delta}}

\def\rrho{\boldsymbol{\rho}}
\newtheorem{lemma}{\hspace*{\parindent}Lemma}
\newtheorem{theorem}{\hspace*{\parindent}Theorem}

\title{Further applications of the $G$ function integral method}
\author{M.A.C. Candezano$^{\rm a}$,  D.B.\:Karp$^{\rm b,c}$\footnote{Corresponding author. E-mail: D. Karp -- \emph{dimkrp@gmail.com},
M.A.C.\:Candezano -- \emph{miguelcaro@mail.uniatlantico.edu.co}, E.\:Prilepkina --  \emph{pril-elena@yandex.ru}}~~and
E.G.\:Prilepkina$^{\rm b,c}$
\\[10pt]
\\
\small{\textit{$\phantom{1}^a$Universidad del Atl\'{a}ntico, Barranquilla, Colombia}}
\\
\small{\textit{$\phantom{1}^b$Far Eastern Federal University, 8
Sukhanova street, Vladivostok, 690950, Russia}}
\\
\small{\textit{$\phantom{1}^c$Institute of Applied Mathematics,
FEBRAS, 7 Radio Street, Vladivostok,  690041, Russia}}}
\date{}

\begin{document}

\maketitle

\begin{abstract}
In our recent work we proposed a generalization of the beta integral method for derivation of the hypergeometric identities which can by analogy be termed ''the $G$ function integral method''.  In this paper we apply this technique to the cubic and the degenerate Miller-Paris transformations to get several new transformation and summation formulas for the generalized hypergeometric functions at a fixed argument.  We further present an alternative approach for reducing the right hand sides resulting from our method to a single hypergeometric function which does not require the use of summation formulas.
\end{abstract}

\bigskip


\bigskip


\bigskip

\section{Introduction and preliminaries}

In our recent paper \cite{KPGFmethod} we proposed a generalization of the beta integral method \cite{KrRao2003} for deriving transformation formulas for hypergeometric functions at a fixed argument. It is based on the following simple idea: the beta density is replaced by a density expressed in terms of Meijer-N{\o}rlund's function $G^{p,0}_{p,p}$ and the Gauss summation theorem for ${}_2F_1$ is replaced by a summation theorem for ${}_{p+1}F_{p}(1)$ with $p\ge2$.  Here  ${}_{p+1}F_{p}$ stands for the generalized hypergeometric function \cite[(2.1.2)]{AAR} and  $G^{p,0}_{p,p}$ is defined in \eqref{eq:G-defined} below.  It is convenient to introduce an extended definition of the  generalized hypergeometric function by
\begin{equation}\label{eq:pFqdefined}
F\left.\!\left(\!\begin{array}{l:}\a \\\b \end{array}\:P\:\right\vert x\right)
=\sum\limits_{n=0}^{\infty}\frac{(a_1)_n(a_2)_n\cdots(a_{p})_n}{(b_1)_n(b_2)_n\cdots(b_q)_nn!}P(n)x^n,
\end{equation}
where $\a=\{a_1,\ldots,a_p\}$, $\b=\{b_1,\ldots,b_q\}$ are complex parameter vectors, $(a)_n=\Gamma(a+n)/\Gamma(a)$ denotes the rising factorial. We found it convenient to omit the indices of the hypergeometric functions, as the dimensions of the parameter vectors are usually clear from the context. However, we will use the traditional notation ${}_pF_{q}$ when dealing with specific numerical values of $p$ and $q$ to make the formulas more accessible to a reader not interested in further details. We will further assume throughout the paper that $b_j$ does not equal a non-positive integer for all $j\in\{1,\ldots,q\}$. The function $P(n)$ in this paper will always be a polynomial of a fixed degree $m$. It is then straightforward to check that
\begin{equation}\label{eq:Pfactorized}
P(n)=P(0)\frac{(1-\llambda)_{n}}{(-\llambda)_{n}},
\end{equation}
where $\llambda=(\lambda_1,\ldots,\lambda_m)$ is the vector of zeros of the polynomial $P$ and the shorthand notation for the product $(-\llambda)_{n}=(-\lambda_1)_n(-\lambda_2)_n\cdots(-\lambda_m)_n$ has been used.  Hence,
\begin{equation}\label{eq:FP-roots}
F\left.\!\left(\!\begin{array}{l:}\a \\\b \end{array}\:P\:\right\vert x\right)
=P(0){}_{p+m}F_{q+m}\left.\!\left(\!\begin{array}{l}\a, 1-\llambda \\\b, -\llambda \end{array}\right\vert x\right)
\end{equation}
- a generalized hypergeometric functions with $m$ unit shifts in parameters.  This extended definition has been recently used by Maier \cite{Maier2019} and is equivalent to the concept of ''hypergeometrization'' introduced  a bit earlier by Blaschke \cite{Blaschke}. We will use both ways of writing $F$ interchangeably.
Omitted argument of the generalized hypergeometric function will signify the unit argument throughout the paper.

The standard symbols $\N$, $\Z$, and $\C$ will be used to denote the sets of natural, integer, and complex numbers, respectively. Similarly to the beta integral method,  our approach starts with a transformation formula of the form
\begin{equation}\label{eq:generaltrans}
F\left.\!\left(\begin{matrix}\aalpha\\\bbeta\end{matrix}\:\right\vert
Mx^{w}\right)
=(1-x)^{\lambda}F\left.\!\left(\begin{matrix}\ddelta\\\ggamma\end{matrix}\:\right\vert
\frac{Dx^{u}}{(1-x)^v}\right)
\end{equation}
valid for $0<x<1$.  Here $\ddelta$, $\ggamma$ and $\lambda$ are functions of $\aalpha$, $\bbeta$; $w,u\in\N$, $v\in\Z$, $M,D$ are constants.  Multiplying this formula by the Meijer-N{\o}rlund function $G^{p,0}_{p,p}$,
defined by the Mellin-Barnes integral of the form
\begin{equation}\label{eq:G-defined}
G^{p,0}_{p,p}\!\left(\!z~\vline\begin{array}{l}\b\\\a\end{array}\!\!\right)\!\!:=
\\
\frac{1}{2\pi{i}}
\int\limits_{\mathcal{L}}\!\!\frac{\Gamma(\a\!+\!s)}{\Gamma(\b+\!s)}z^{-s}ds,
\end{equation}
and integrating term-wise we established the ''master'' lemma  below \cite[Lemma~1]{KPGFmethod}.
Details regarding the choice of the contour ${\mathcal{L}}$ can be found in many standard reference books \cite[section~5.2]{LukeBook}, \cite[16.17]{NIST}, \cite[8.2]{PBM3} and our papers \cite{KLJAT2017,KPSIGMA}, which also contain a list of properties of $G^{p,0}_{p,p}$.
\begin{lemma}\label{lm:master1}
Assume that \eqref{eq:generaltrans} holds for  $x\in(0,1)$. Suppose further that $\ddelta$ or $\a$ contains a negative integer or $v=0$, $D=1$, and
\begin{equation}\label{eq:v0condition}
\Re(\a)>0~\&~\Re(s(\a,\b)+\lambda)>0~\&~\Re(s(\a,\b)+s(\ggamma,\ddelta)+\lambda)>0,
\end{equation}
where $s(\a,\b)=\sum_{j=1}^{p}(b_{j}-a_{j})$ is the parametric excess. Then
\begin{equation}\label{eq:master1}
F\left.\!\!\left(\!\begin{matrix}\aalpha,\Delta(\a,w)\\\bbeta,\Delta(\b,w)\end{matrix}\right\vert M\right)
=\sum\limits_{k=0}^{\infty}\frac{(\ddelta)_k(\a)_{uk}D^k}{(\ggamma)_k(\b)_{uk}k!}
F\left.\!\!\left(\!\begin{matrix}-\lambda+vk,\a+uk\\\b+uk\end{matrix}\right.\right),
\end{equation}
where $\Delta(a,w)=(a/w,a/w+1/w,\ldots,a/w+(w-1)/w)$.
\end{lemma}

In \cite{KPGFmethod} we applied our method to a number of transformations with $w,u,v\in\{-1,0,1,2\}$ including Euler-Pfaff, Miller-Paris and many quadratic transformations.    The purpose of this note is threefold. First, we apply the method to the cubic and the degenerate Miller-Paris transformations; second, we propose an alternative way to handle the expression on the right hand side of \eqref{eq:master1}; finally, we will show how transformation formulas obtained by $G$ function integral method can be used to derive summation formulas for the generalized hypergeometric functions including those with with non-linearly constrained parameters.

Before moving forward to these topics let us cite Lemma~2 from \cite{KPGFmethod}, whose particular cases will be used extensively to sum the hypergeometric function on the right hand side of \eqref{eq:master1}.
\begin{lemma}\label{lm:IPDsummation}
Suppose $l\in\N$, $\hh=(h_1,\ldots,h_l)\in\C^{l}$, $\pp=(p_1,\ldots,p_l)\in\N^{l}$, $p=p_1+\cdots+p_l$, $u\in\N$, $v\in\Z$.
Then for $k\in\N$ such that $\Re(e+\lambda-d-p-vk)>0$ or if hypergeometric function $F$ terminates, we have
\begin{equation}\label{eq:IPDsummation}
F\!\left(\begin{matrix}-\lambda+vk,d+uk,\hh+\pp+uk\\e+uk,\hh+uk\end{matrix}\right)
=\frac{(-1)^{vk}\Gamma(e+\lambda-d)\Gamma(1+d-e-\lambda)\Gamma(e+uk)Y_p(u,v;k)}{(\hh+uk)_{\pp}\Gamma((u-v)k+e+\lambda)
\Gamma(vk+d-e-\lambda+p+1)},
\end{equation}
where $(\hh+uk)_{\pp}=(h_1+uk)_{p_1}\cdots(h_l+uk)_{p_l}$ and
\begin{equation}\label{eq:Ypolynomial}
Y_p(u,v;t)=\frac{(\hh-d)_{\pp}}{\Gamma(e-d)}\sum_{j=0}^{p}\frac{(d-e+1)_{j}}{j!}
F\!\left(\begin{matrix}-j,1-\hh+d\\1-\hh+d-\pp\end{matrix}\right)(ut+d)_{j}(vt+d-e-\lambda+j+1)_{p-j}
\end{equation}
is a polynomial of degree $p$.
\end{lemma}

\noindent\textbf{Remark.} If $p=l=1$ the polynomial $Y_p(u,v;t)$ reduces to
\begin{equation}\label{eq:Y1}
Y_1(u,v;t)=\frac{v(h-d)-u(e-d-1)}{\Gamma(e-d)}t-\frac{(h-d)\lambda+(e-d-1)h}{\Gamma(e-d)}
\end{equation}
with the root
\begin{equation}\label{eq:Y1root}
\xi=\frac{(h-d)\lambda+(e-d-1)h}{(h-d)v-(e-d-1)u}.
\end{equation}

Denote $\Delta(z,m)_{k}=(z/m)_{k}((z+1)/m)_{k}\cdots((z+m-1)/m)_{k}$. We will need several particular cases of the above lemma which are easily derived from \eqref{eq:IPDsummation} using the identities
$$
\Gamma(z-n)=\frac{(-1)^n\Gamma(z)}{(1-z)_{n}},~~~(z)_{2k}=4^k\Delta(z,2)_{k},~~~(z)_{3k}=27^k\Delta(z,3)_{k}.
$$
For $(u,v)=(1,3)$ we have
\begin{equation}\label{eq:IPD1-13}
F\!\left(\begin{matrix}-\lambda+3k,d+k,\hh+\pp+k\\e+k,\hh+k\end{matrix}\right)
=\frac{\Gamma(e+\lambda-d)\Gamma(e)(e)_k(\hh)_{k}\Delta(1-e-\lambda,2)_k(-4/27)^kY_p(1,3;k)}{(\hh)_{\pp}\Gamma(e+\lambda)(1+d-e-\lambda)_{p}(\hh+\pp)_{k}\Delta(1+d-e-\lambda+p,3)_{k}}.
\end{equation}
The case $(u,v)=(2,3)$ takes the form:
\begin{multline}\label{eq:IPD1-23}
F\!\left(\begin{matrix}-\lambda+3k,d+2k,\hh+\pp+2k\\e+2k,\hh+2k\end{matrix}\right)
\\
=\frac{\Gamma(e+\lambda-d)\Gamma(e)\Delta(e,2)_k\Delta(\hh,2)_{k}(1-e-\lambda)_{k}(4/27)^kY_p(2,3;k)}{(\hh)_{\pp}\Gamma(e+\lambda)(1+d-e-\lambda)_{p}\Delta(\hh+\pp,2)_{k}\Delta(1+d-e-\lambda+p,3)_{k}}.
\end{multline}
The case $(u,v)=(3,2)$ is given by
\begin{multline}\label{eq:IPD1-32}
F\!\left(\begin{matrix}-\lambda+2k,d+3k,\hh+\pp+3k\\e+3k,\hh+3k\end{matrix}\right)
\\
=\frac{\Gamma(e+\lambda-d)\Gamma(e)\Delta(e,3)_k\Delta(\hh,3)_{k}(27/4)^kY_p(3,2;k)}{(\hh)_{\pp}\Gamma(e+\lambda)(1+d-e-\lambda)_{p}\Delta(\hh+\pp,3)_{k}(e+\lambda)_{k}\Delta(1+d-e-\lambda+p,2)_{k}}.
\end{multline}
Finally, for $(u,v)=(1,-2)$  we obtain:
\begin{multline}\label{eq:IPD1-1-2}
F\!\left(\begin{matrix}-\lambda-2k,d+k,\hh+\pp+k\\e+k,\hh+k\end{matrix}\right)
\\
=\frac{\Gamma(e+\lambda-d)\Gamma(e)(e)_k(\hh)_{k}(4/27)^k\Delta(e+\lambda-d-p,2)_{k}Y_p(1,-2;k)}{(\hh)_{\pp}\Gamma(e+\lambda)(1+d-e-\lambda)_{p}(\hh+\pp)_{k}\Delta(e+\lambda,3)_{k}}.
\end{multline}

\section{Cubic transformations}

The following lemma based on the Gessel-Stanton identity \cite[(1.9)]{GesselStant} yields one more summation formula for the case  $(u,v)=(1,3)$.
\begin{lemma}\label{lm:GSsum}
For any $n,k\in\N$, $k\le{n}$, we have
\begin{multline}\label{eq:GSsum}
{}_{3}F_{2}\left.\!\left(\begin{matrix}-\lambda+3k,-(\lambda+1)/3+k,-n+k\\
(5-\lambda)/3+k,2n+2-\lambda+k\end{matrix}\right.\right)
\\
=\frac{(2-\lambda)(n+1)(3n+2)!(-1/2-n)_{k}(2-\lambda+2n)_{k}((5-\lambda)/3)_{k}(-4/27)^{k}}
{(2-\lambda+3n)(2n+2)!(2-\lambda+2n)_{n}(-2/3-n)_{k}(-1/3-n)_k((2-\lambda)/3)_{k}}.
\end{multline}
\end{lemma}
\textbf{Proof.} According to $s=-3$ case of \cite[(1.9)]{GesselStant} (see also \cite[(8.12)]{Koepf}), we have:
$$
{}_{3}F_{2}\left.\!\!\left(\begin{matrix}-\lambda+3k,-(\lambda+1)/3+k,-n+k\\
(5-\lambda)/3+k,2n+2-\lambda+k\end{matrix}\right.\right)
=\frac{(-2-3(n-k))_{n-k}((2-\lambda)/3+k)(n-k+1)}{(-1+\lambda-3n))_{n-k}(2-\lambda)/3+n)}.
$$
Next, apply the easily verifiable identities
$$
(-2-3n+3k)_{n-k}=(-1)^{n}\frac{(3n+2)!(-1-n)_{k}(-1/2-n)_{k}4^k}{(2n+2)!(-2/3-n)_{k}(-1/3-n)_{k}(-n)_{k}27^k},
$$
$$
(-1+\lambda-3n)_{n-k}=(-1)^{n+k}\frac{(2-\lambda+2n)_{n}}{(2-\lambda+2n)_{k}}
$$
to get \eqref{eq:GSsum}.$\hfill\square$

Combining Lemma~\ref{lm:master1} with Lemma~\ref{lm:GSsum} and summation formulas \eqref{eq:IPD1-13}-\eqref{eq:IPD1-1-2} we obtain a number of transformation formulas for terminating generalized hypergeometric functions, none of which could be immediately located in the literature.  We will present each formula in a separate theorem.  Recall that $\Delta(z,m)=(z/m,z/m+1/m,\ldots,z/m+(m-1)/m)$ and bottom parameters are always assumed to satisfy the restriction of not being equal to non-positive integers.
\begin{theorem}\label{th:Bailey-GS}
For $n\in\N$ we have
\begin{multline}\label{eq:Bailey-GS}
{}_{5}F_{4}\left.\!\!\left(\begin{matrix}\alpha,2\beta-\alpha-1,\alpha+2-2\beta,(\alpha-1)/3,-n\\
\beta,\alpha-\beta+3/2,(\alpha+5)/3,2n+\alpha+2\end{matrix}\right|\frac{1}{4}\right)
\\
=\frac{(\alpha+2)\Gamma(3n+4)\Gamma(2n+\alpha+2)}{3\Gamma(2n+3)\Gamma(3n+\alpha+3)}
{}_{5}F_{4}\left.\!\!\left(\begin{matrix}\Delta(\alpha-1,3),-1/2-n,-n\\
\beta,\alpha-\beta+3/2,-1/3-n,-2/3-n\end{matrix}\right.\right).
\end{multline}
\end{theorem}
\textbf{Proof.} According to \cite[p.185]{AAR} the following cubic transformation due to Bailey holds true for $0<x<1$:
\begin{equation}\label{eq:BaileyCubic}
{}_{3}F_{2}\left.\!\left(\begin{matrix}\alpha, 2\beta-\alpha-1,\alpha+2-2\beta\\
\beta,\alpha+3/2-\beta\end{matrix}\right|\frac{x}{4}\right)
=(1-x)^{-\alpha}{}_{3}F_{2}\left.\!\left(\begin{matrix}\alpha/3, (\alpha+1)/3,(\alpha+2)/3\\
\beta,\alpha+3/2-\beta\end{matrix}\right|\frac{-27x}{4(1-x)^3}\right).
\end{equation}
Then we can apply Lemma~\ref{lm:master1} with $M=1/4$, $u=1$, $v=3$, $\lambda=-\alpha$ and $D=-27/4$.  This yields
\begin{equation}\label{eq:Bailey-GSproof}
F\left.\!\!\left(\!\begin{matrix}\alpha, 2\beta-\alpha-1,\alpha+2-2\beta,\a\\\beta,\alpha+3/2-\beta,\b\end{matrix}\right\vert \frac{1}{4}\right)
=\sum\limits_{k=0}^{\infty}\frac{\Delta(\alpha,3)_{k}(\a)_{k}(-27/4)^k}{(\beta)_{k}(\alpha+3/2-\beta)_{k}(\b)_{k}k!}
F\left.\!\!\left(\!\begin{matrix}\alpha+3k,\a+k\\\b+k\end{matrix}\right.\right).
\end{equation}
By choosing $a_1=(1-\alpha)/3$, $a_2=-n$, $n\in\N$, $b_1=(5+\alpha)/3$, $b_2=2n+2+\alpha$  we are in the position to apply Lemma~\ref{lm:GSsum} which, after some cancelations, leads to  \eqref{eq:Bailey-GS}. $\hfill\square$

\begin{theorem}\label{th:Bailey-IPD}
For $-d\in\N$ we have
\begin{multline}\label{eq:Bailey-IPD}
F\left.\!\!\left(\begin{matrix}\alpha,2\beta-\alpha-1,\alpha+2-2\beta,d,\hh+\pp\\
\beta,\alpha-\beta+3/2,e,\hh\end{matrix}\right|\frac{1}{4}\right)
\\
=\frac{\Gamma(e-\alpha-d)\Gamma(e)}{(\hh)_{\pp}\Gamma(e-\alpha)(1+d+\alpha-e)_{p}}
 F\left.\!\left(\!\begin{array}{l:}\Delta(\alpha,3),\Delta(1+\alpha-e,2),d
\\
\beta,\alpha-\beta+3/2,\Delta(1+d+\alpha-e+p,3)\end{array}\:Y_p(1,3)\:\right.\!\!\right),
\end{multline}
where $Y_p(1,3)$ is defined in \eqref{eq:Ypolynomial} with $\lambda=-\alpha$.
\end{theorem}
\textbf{Proof.}  Follow the proof of Theorem~\ref{th:Bailey-GS} up to formula \eqref{eq:Bailey-GSproof}. Then choose $a_1=d$, $\a_{[1]}=\hh+\pp$,
$b_1=e$, $\b_{[1]}=\hh$ and apply formula \eqref{eq:IPD1-13}.  Here $\a_{[1]}$ is a shorthand notation for the vector $\a$ with the first component removed. $\hfill\square$

\begin{theorem}\label{th:GScubic-GS}
For $n\in\N$ we have
\begin{multline}\label{eq:GScubic-GS}
{}_{4}F_{3}\left.\!\!\left(\begin{matrix}1,(\alpha-1)/3,\alpha-2\beta,-n\\
\beta+1,(\alpha+5)/3,2n+\alpha+2\end{matrix}\right.\right)
\\
=\frac{(\alpha+2)\Gamma(3n+4)\Gamma(2n+\alpha+2)}{3\Gamma(2n+3)\Gamma(3n+\alpha+3)}
{}_{6}F_{5}\left.\!\!\left(\begin{matrix}\Delta(\alpha-1,3),\beta,-1/2-n,-n\\
\Delta(\alpha,2),\beta+1,-1/3-n,-2/3-n\end{matrix}\right.\right).
\end{multline}
The function ${}_{6}F_{5}$ on the right hand side is Saalsh\"{u}tzian \emph{(}or $1$-balanced\emph{)}.
\end{theorem}
\textbf{Proof.} According to \cite[(5.13)]{GesselStant} we have:
\begin{equation}\label{eq:GScubic}
{}_{2}F_{1}\left.\!\!\!\left(\begin{matrix}1,\alpha-2\beta\\
\beta+1\end{matrix}\right|x\right)
=(1-x)^{-\alpha}{}_{4}F_{3}\left.\!\!\!\left(\begin{matrix}\alpha/3, (\alpha+1)/3,(\alpha+2)/3,\beta\\
\alpha/2,(\alpha+1)/2,\beta+1\end{matrix}\right|\frac{-27x}{4(1-x)^3}\right).
\end{equation}
Then we can apply Lemma~\ref{lm:master1} with $M=1$, $u=1$, $v=3$, $\lambda=-\alpha$ and $D=-27/4$.  Then setting
 $a_1=(1-\alpha)/3$, $a_2=-n$, $n\in\N$, $b_1=(5+\alpha)/3$, $b_2=2n+2+\alpha$  we can apply Lemma~\ref{lm:GSsum} to sum the hypergeometric function of the right hand side. This leads immediately to \eqref{eq:GScubic-GS}.$\hfill\square$

\begin{theorem}\label{th:GScubic-IPD}
For $-d\in\N$ we have
\begin{multline}\label{eq:GScubic-IPD}
F\left.\!\!\left(\begin{matrix}1,\alpha-2\beta,d,\hh+\pp\\
\beta+1,e,\hh\end{matrix}\right.\right)
\\
=\frac{\Gamma(e-\alpha-d)\Gamma(e)}{(\hh)_{\pp}\Gamma(e-\alpha)(1+d+\alpha-e)_{p}}
 F\left.\!\left(\!\begin{array}{l:}\Delta(\alpha,3),\beta,\Delta(1+\alpha-e,2),d
\\
\Delta(\alpha,2),\beta+1,\Delta(1+d+\alpha-e+p,3)\end{array}\:Y_p(1,3)\:\right.\!\!\right),
\end{multline}
where $Y_p(1,3)$ is defined in \eqref{eq:Ypolynomial} with $\lambda=\alpha$.
\end{theorem}
\textbf{Proof.} Use \eqref{eq:GScubic} in Lemma~\ref{lm:master1} similarly to the proof of Theorem~\ref{th:GScubic-GS}. Then apply  formula \eqref{eq:IPD1-13} to sum the hypergeometric function on the right hand side.$\hfill\square$

\begin{theorem}\label{th:AARcubic-IPD}
For $-d\in\N$ we have
\begin{multline}\label{eq:AARcubic-IPD}
F\left.\!\!\left(\begin{matrix}\alpha,\beta-1/2,\alpha-\beta+1,d,\hh+\pp\\
2\beta,2\alpha-2\beta+2,e,\hh\end{matrix}\right|4\right)
=\frac{\Gamma(e-\alpha-d)\Gamma(e)}{(\hh)_{\pp}\Gamma(e-\alpha)(1+d+\alpha-e)_{p}}
\\
 \times\! F\left.\!\!\left(\!\!\begin{array}{l:}\Delta(\alpha,3),\Delta(e,2),\Delta(\hh,2),1+\alpha-e,d,\hh+\pp
\\
\beta,\alpha-\beta+3/2,e,\hh,\Delta(\hh+\pp,2),\Delta(1+d+\alpha-e+p,3)\end{array}\:Y_p(2,3)\:\right.\!\!\right),
\end{multline}
where $Y_p(2,3)$ is defined in \eqref{eq:Ypolynomial} with $\lambda=\alpha$.
\end{theorem}
\textbf{Proof.} According to \cite[p.185]{AAR} the following cubic transformation due to Bailey holds true for $0<x<1$:
\begin{equation}\label{eq:AAR185}
{}_{3}F_{2}\left.\!\left(\begin{matrix}\alpha, \beta-1/2,\alpha+1-\beta\\
2\beta,2\alpha+2-2\beta\end{matrix}\right|4x\right)
=(1-x)^{-\alpha}{}_{3}F_{2}\left.\!\left(\begin{matrix}\Delta(\alpha,3)\\
\beta,\alpha+3/2-\beta\end{matrix}\right|\frac{27x^2}{4(1-x)^3}\right).
\end{equation}
Application of Lemma~\ref{lm:master1} and  formula \eqref{eq:IPD1-23} completes the proof.$\hfill\square$

\textbf{Remark.} Bailey's cubic transformation \eqref{eq:AAR185} has been recently extended by Maier in \cite[Theorems~3.3, 3.6, 3.9]{Maier2019}.  These extensions can be used in place of  \eqref{eq:AAR185} to derive generalizations of \eqref{eq:AARcubic-IPD}.

\begin{theorem}\label{th:GS5.20-IPD}
For $-d\in\N$ we have
\begin{multline}\label{eq::GS5.20-IPD}
F\left.\!\!\left(\begin{matrix}3\alpha,3\alpha+1/2,d,\hh+\pp\\
6\alpha+1,e,\hh\end{matrix}\right|\frac{4}{3}\right)
=\frac{\Gamma(e-2\alpha-d)\Gamma(e)}{(\hh)_{\pp}\Gamma(e-2\alpha)(1+d+2\alpha-e)_{p}}
\\
 \times\! F\left.\!\!\left(\!\!\begin{array}{l:}\alpha,\alpha+1/2,\Delta(d,3)
\\
2\alpha+1,e-2\alpha,\Delta(1+d+2\alpha-e+p,2)\end{array}\:Y_p(3,2)\:\right.\!\!\right),
\end{multline}
where $Y_p(3,2)$ is defined in \eqref{eq:Ypolynomial} with $\lambda=2\alpha$.
\end{theorem}
\textbf{Proof.}
According to \cite[(5.20)]{GesselStant} we have:
$$
{}_{2}F_{1}\left.\!\!\!\left(\begin{matrix}3\alpha,3\alpha+1/2\\
6\alpha+1\end{matrix}\right|\frac{4x}{3}\right)
=(1-x)^{-2\alpha}
{}_{2}F_{1}\left.\!\!\!\left(\begin{matrix}\alpha, \alpha+1/2\\
2\alpha+1\end{matrix}\right|\frac{4x^3}{27(1-x)^2}\right).
$$
Application of Lemma~\ref{lm:master1} and  formula \eqref{eq:IPD1-32} completes the proof.$\hfill\square$

For $p=1$ formula \eqref{eq::GS5.20-IPD} takes the form:
\begin{multline*}
{}_4F_{3}\left.\!\!\left(\begin{matrix}3\alpha,3\alpha+1/2,d,h+1\\
6\alpha+1,e,h\end{matrix}\right|\frac{4}{3}\right)
=\frac{\Gamma(e-2\alpha-d)\Gamma(e)(2\alpha(h-d)-(e-d-1)h)}{h\Gamma(e-d)\Gamma(e-2\alpha)(1+d+2\alpha-e)}
\\
 \times{}_6F_{5}\left.\!\left(\begin{matrix}\alpha,\alpha+1/2,\Delta(d,3),\xi+1
\\
2\alpha+1,e-2\alpha,\Delta(2+d+2\alpha-e,2),\xi\end{matrix}\right.\!\right),
\end{multline*}
where
$$
\xi=\frac{2\alpha(h-d)-(e-d-1)h}{2(h-d)-3(e-d-1)}.
$$

\begin{theorem}\label{th:GS5.18-IPD}
Suppose $\alpha$ or $d$ is a negative integer. Then
\begin{equation}\label{eq:GS5.18-IPD}
F\left.\!\!\left(\begin{matrix}3\alpha,-3\alpha,d,\hh+\pp\\
1/2,e,\hh\end{matrix}\right|\frac{3}{4}\right)
=\frac{\Gamma(e-d)}{(\hh)_{\pp}(1+d-e)_{p}}
F\left.\!\!\left(\!\begin{array}{l:}\alpha,-\alpha,d,\Delta(e-d-p,2)
\\
1/2,\Delta(e,3)\end{array}\:Y_p(1,-2)\:\right.\!\!\right),
\end{equation}
where $Y_p(1,-2)$ is defined in \eqref{eq:Ypolynomial} with $\lambda=0$.
\end{theorem}
\textbf{Proof.} We start with the transformation \cite[(5.18)]{GesselStant}
$$
{}_{2}F_{1}\left.\!\!\!\left(\begin{matrix}3\alpha,-3\alpha\\
1/2\end{matrix}\right|\frac{3x}{4}\right)
={}_{2}F_{1}\left.\!\!\!\left(\begin{matrix}\alpha, -\alpha\\
1/2\end{matrix}\right|\frac{27x(1-x)^2}{4}\right)
$$
playing the role of \eqref{eq:generaltrans} in Lemma~\ref{lm:master1}.  Then use formula \eqref{eq:IPD1-1-2} to sum the hypergeometric function on the right hand side.$\hfill\square$

In view of \eqref{eq:Y1} and \eqref{eq:Y1root}, transformation \eqref{eq:GS5.18-IPD} takes a particularly simple form for $p=1$:
$$
{}_4F_{3}\left.\!\!\left(\begin{matrix}3\alpha,-3\alpha,d,h+1\\
1/2,e,h\end{matrix}\right|\frac{3}{4}\right)
={}_6F_{5}\left.\!\!\left(\begin{matrix}\alpha,-\alpha,d,\Delta(e-d-1,2),\xi+1
\\
1/2,\Delta(e,3),\xi\end{matrix}\right.\!\right),
$$
where $\xi=(e-d-1)h/(2h+e-3d-1)$ is the negated root of $Y_1(1,-2)$.
Further, setting $e=d+1+\varepsilon$ and letting $\varepsilon\to0$ after some algebra we arrive at ($\alpha\in\N$):
$$
{}_4F_{3}\left.\!\!\left(\begin{matrix}3\alpha,-3\alpha,d,h+1\\
1/2,d+1,h\end{matrix}\right|\frac{3}{4}\right)
=\frac{d}{h}+\frac{h-d}{h}{}_4F_{3}\left.\!\!\left(\begin{matrix}\alpha,-\alpha,d,1
\\
\Delta(d+1,3)\end{matrix}\right.\!\right).
$$
Note that ${}_4F_{3}$ on the right hand side is Saalsch\"{u}tzian (i.e. $1$-balanced) while ${}_4F_{3}$ on the left hand side is $1/2$-balanced.

Bailey's cubic transformations have been recently extended by Maier in \cite{Maier2019}. These extensions can be combined with Lemma~\ref{lm:IPDsummation} to get generalizations of Theorems~\ref{th:Bailey-IPD} and  \ref{th:AARcubic-IPD}.  Three Maier's transformations \cite[Theorems~3.2, 3.5, 3.8]{Maier2019}  can also be combined with Lemma~\ref{lm:GSsum}.  We will restrict our attention to a combination of Lemma~\ref{lm:GSsum} with the transformation \cite[Theorems~3.2]{Maier2019}
\begin{multline}\label{eq:MaierTh3.2}
F\left.\!\left(\!\begin{array}{l:}\alpha,1/2-r-\beta,1/2-r+\beta
\\
(\alpha+\beta+r)/2+3/4, (\alpha-\beta+r)/2+3/4\end{array}\:Q^{(3)}_{r}\:\right|\frac{x}{4}\right)
\\
=(1-x)^{-\alpha}{}_{3}F_{2}\left.\!\left(\begin{matrix}\Delta(\alpha,3)\\
(\alpha+\beta+r)/2+3/4, (\alpha-\beta+r)/2+3/4\end{matrix}\right|\frac{-27x}{4(1-x)^3}\right).
\end{multline}
Here the $2r$-degree polynomial $Q^{(3)}_{r}$ is given by
$$
Q^{(3)}_{r}(x)=\frac{((\frac{1}{2}-r+\beta-x)/2)_{r}((\frac{1}{2}-r-\beta-x)/2)_{r}}{4^{-r}(1/2+\beta)_{r}(1/2-\beta)_{r}}
{}_{3}F_{2}\left.\!\left(\begin{matrix}-x,(x+\alpha)/2,-r\\
(\frac{1}{2}-r+\beta-x)/2,(\frac{1}{2}-r-\beta-x)/2\end{matrix}\right.\right).
$$
This leads to the following generalization of Theorem~\ref{th:Bailey-GS}.
\begin{theorem}\label{th:Maier3.2-GS}
For $n\in\N$ we have
\begin{multline}\label{eq:Maier3.2-GS}
F\left.\!\!\left(\begin{array}{l:}\alpha,1/2-r-\beta,1/2-r+\beta,(\alpha-1)/3,-n\\
(\alpha+\beta+r)/2+3/4, (\alpha-\beta+r)/2+3/4,(\alpha+5)/3,2n+\alpha+2\end{array}\:Q^{(3)}_{r}\:\right|\frac{1}{4}\right)
\\
=\frac{(\alpha+2)\Gamma(3n+4)\Gamma(2n+\alpha+2)}{3\Gamma(2n+3)\Gamma(3n+\alpha+3)}
{}_{5}F_{4}\left.\!\!\left(\begin{matrix}\Delta(\alpha-1,3),-1/2-n,-n\\
(\alpha+\beta+r)/2+3/4, (\alpha-\beta+r)/2+3/4,-1/3-n,-2/3-n\end{matrix}\right.\right).
\end{multline}
\end{theorem}
\textbf{Proof.} Apply Lemma~\ref{lm:master1} to transformation \eqref{eq:MaierTh3.2} and use Lemma~\ref{lm:GSsum} to sum the hypergeometric functions on the right hand side.$\hfill\square$

\section{Degenerate Miller-Paris transformations}

Miller-Paris transformations are extensions of Euler's transformations for the Gauss hypergeometric function \cite[Theorem~2.2.5]{AAR} to generalized hypergeometric functions of higher-order having integral parameter differences (IPD-type).  They were developed in a series of papers published over last 15 years, the most general form was presented in a seminal paper \cite{MP2013} by Miller and Paris.  In our recent articles \cite{KPResults2019,KarpPril1} we extended these transformations to the previously prohibited valued of parameters and gave denomination ''degenerate Miller-Paris transformations'' to the resulting identities. In this section we apply the $G$ function integral method to some degenerate Miller-Paris transformations.  As these transformations are not of the form \eqref{eq:generaltrans}, we cannot use Lemma~\ref{lm:master1}, so we will follow the method explicitly.  As we mentioned earlier the essence of the method is to multiply a known transformation by the Meijer-N{\o}rlund function $G^{p,0}_{p,p}$ function defined in \eqref{eq:G-defined} and integrate it from $0$ to $1$.  To perform the term-wise integration we will need the integral evaluation \cite[p.50]{KLJAT2017}
\begin{equation}\label{eq:Gpowerint}
\frac{\Gamma(\b)}{\Gamma(\a)}\int_0^1x^{\nu}(1-x)^{\mu}G^{p,0}_{p,p}\!\left(\!x\left|\!\begin{array}{l}\b-1\\\a-1\end{array}\right.\!\!\!\right)dx
=\frac{(\a)_{\nu}}{(\b)_{\nu}}{}_{p+1}F_{p}\left.\!\!\left(\!\begin{matrix}-\mu,\a+\nu\\\b+\nu\end{matrix}\right.\right),
\end{equation}
where for arbitrary $\nu$ the Pochhammer's symbol is defined by $(a)_{\nu}=\Gamma(a+\nu)/\Gamma(a)$. The above formula is true if $\Re(\a+\nu)>0$ and $\Re(s(\a,\b)+\mu)>0$ (recall that $s(\a,\b)=\sum_{j=1}^{p}(b_j-a_j)$).
We now apply this technique to the degenerate Miller-Paris transformation found in \cite{KarpPril1}.
Define  $\m=(m_1,\ldots,m_r)\in\N^r$, $m=m_1+m_2+\ldots+m_r$ and  $\f=(f_1,\ldots,f_r)\in\C^{r}$.  
We will reserve the symbols $\f$ and $\m$ for the degenerate Miller-Paris transformations throughout the rest of the paper.

\begin{theorem}\label{th:Degenerate1}
Suppose $\Re(e-a-d-p-m+1)>0$, $\Re(e-d-p)>0$. Then
\begin{multline}\label{eq:Kartransfor1}
F\left.\!\left(\begin{matrix}a,b,d,\f+\m,\hh+\pp\\
b+1,e,\f,\hh\end{matrix}\right.\right)
=\Omega\cdot F\left.\!\left(\!\begin{array}{l:}1,
b+1-a,d \\b+1,e-a+1 \end{array}\:Y_p(1,0)\:\right.\right)
\\
+\sum\limits_{l=0}^{m-1}\beta_l\frac{(a)_l(d)_l(\hh+\pp)_l}{(e)_l(\hh)_l}
F\left.\!\left(\begin{matrix}a+l,d+l,\hh+\pp+l\\e+l,\hh+l\end{matrix}\:\right.\right),
\end{multline}
where the polynomial $Y_p(1,0)$ is defined in \eqref{eq:Ypolynomial} with $\lambda=1-a$ and
\begin{equation}\label{eq:Bl1}
\Omega\!=\!\frac{(\f-b)_\m\Gamma(e)\Gamma(e-a-d+1)}{(\f)_\m(\hh)_\pp(d+a-e)_p\Gamma(e+1-a)},
~~\beta_l=\frac{(-1)^{l}}{l!}{}F\!\!\left(\begin{matrix}-l,
b,\f+\m\\ b+1,\f\end{matrix}\right)
-\frac{(-1)^{l}(\f-b)_\m}{(b+1)_l(\f)_\m}.
\end{equation}
\end{theorem}
\textbf{Proof.} According to \cite[Teorem~3]{KarpPril1} \begin{equation}\label{eq:Kar112}
F\left.\!\left(\begin{matrix}a, b,\f+\m\\
b+1,\f\end{matrix}\right\vert x\right)
=\frac{(\f-b)_\m}{(\f)_\m}(1-x)^{1-a}{}_{2}F_{1}\left.\!\!\left(\begin{matrix}1,
b+1-a\\b+1\end{matrix}\right\vert
x\right)+\sum\limits_{l=0}^{m-1}\beta_l(a)_lx^l(1-x)^{-a-l}
\end{equation}
with $\beta_l$  defined in \eqref{eq:Bl1}. Suppose $\hh=(h_1,\ldots,h_l)$ is a complex vector,  $\pp=(p_1,\ldots,p_l)$ comprises non-negative integers, $p=p_1+\ldots+p_l.$  To prove the theorem, multiply both sides of  \eqref{eq:Kar112} by
$$
\frac{\Gamma(e)}{\Gamma(d)(\hh)_\pp} G^{l+1,0}_{l+1,l+1}\!\left(\!x\left|\!\begin{array}{l}e-1,\hh-1\\d-1,\hh+\pp-1\end{array}\right.\!\!\!\right)
$$
and integrate term-wise from $0$ to $1$.  Applying (\ref{eq:Gpowerint}) we then have
\begin{multline}\label{eq:Gint1}
\frac{\Gamma(e)}{\Gamma(d)(\hh)_\pp}\int_0^1 x^l(1-x)^{-a-l}G^{l+1,0}_{l+1,l+1}\!\left(\!x\left|\!\begin{array}{l}e-1,\hh-1\\d-1,\hh+\pp-1\end{array}\right.\!\!\!\right) dx
\\
=\frac{(d)_l(\hh+\pp)_l}{(e)_l(\hh)_l}F\left.\!\left(\begin{matrix}a+l,d+l,\hh+\pp+l\\e+l,\hh+l\end{matrix}\:\right.\right).
\end{multline}
Setting $\aalpha=(a,b,\f+\m)$, $\bbeta=(b+1,\f)$, on the left hand side we obtain:
\begin{equation}\label{eq:Gpowerint1}
\frac{\Gamma(e)}{\Gamma(d)(\hh)_\pp}\int_0^1 G^{l+1,0}_{l+1,l+1}\!\left(\!x\left|\!\begin{array}{l}e-1,\hh-1\\d-1,\hh+\pp-1\end{array}\right.\!\!\!\right) F\left.\!\left(\begin{matrix}\aalpha\\\bbeta\end{matrix}\:\right\vert
x\right)dx=F\left.\!\left(\begin{matrix}\aalpha,d,\hh+\pp\\\bbeta,e,\hh\end{matrix}\:\right.\right).
\end{equation}
Further, term-wise integration and Lemma \ref{lm:IPDsummation} lead to the following evaluation:
\begin{multline}\label{eq:Gpowerint2}
\frac{\Gamma(e)}{\Gamma(d)(\hh)_\pp}\int_0^1 G^{l+1,0}_{l+1,l+1}\!\left(\!x\left|\!\begin{array}{l}e-1,\hh-1\\d-1,\hh+\pp-1\end{array}\right.\!\!\!\right) (1-x)^\lambda{}_2F_{1}\left.\!\left(\begin{matrix}1,
b+1-a\\b+1\end{matrix}\:\right\vert x\right)dx
\\
=\frac{\Gamma(e)\Gamma(e+\lambda-d)}{(\hh)_\pp(d-e-\lambda+1)_p\Gamma(e+\lambda)}
F\left.\!\left(\!\begin{array}{l:}1,b+1-a,d
\\
b+1,e+\lambda \end{array}\:Y_p(1,0)\:\right.\right),
\end{multline}
where the  polynomial $Y_p(1,0)$ is defined by \eqref{eq:Ypolynomial}.~~~~~$\square$

The previous theorem can be further generalized by substituting the bottom parameter $b+1$ by $b+k$ with $k\ge2$, as follows.
\begin{theorem}\label{th:Degenerate2}
Suppose $k\in\N$, $k\ge2$,  $\Re(e-a-d-p-m+1)>0$ and $\Re(e-d-p)>0$. Then
\begin{multline}\label{eq:Kar2transform}
F\left.\!\!\left(\begin{matrix}a,b,d,\f+\m,\hh+\pp\\b+k, e,\f,\hh\end{matrix}\right.\right)
\!=\!\Lambda\cdot
F\left.\!\left(\!\begin{array}{l:}1,b-a+1,d \\b+k,e-a+1 \end{array}\:T_{k-1}\cdot{Y_p(1,0)}\:\right.\right)
\\
+\sum\limits_{q=1}^k\frac{(-1)^{q-1}(b)_{k}}{(b+q-1)(q-1)!(k-q)!}
\sum\limits_{l=0}^{m-1}\sigma_{ql}\frac{(a)_l(d)_l(\hh+\pp)_l}{(e)_l(\hh)_l}
F\left.\!\left(\begin{matrix}a+l,d+l,\hh+\pp+l\\e+l,\hh+l\end{matrix}\:\right.\right),
\end{multline}
where the polynomial $Y_p(1,0)$ is defined by \eqref{eq:Ypolynomial} with $\lambda=1-a$,
$$
\Lambda=\frac{\Gamma(b-a+1)\Gamma(e)\Gamma(e+1-a-d)}{(\f)_\m(\hh)_\pp(d-e+a)_p\Gamma(b)\Gamma(e-a+1)},
$$
$$
T_{k-1}(z)=\sum\limits_{q=1}^k
\frac{(-1)^{q-1}(\f-b-q+1)_\m\Gamma(b+q-1)}{\Gamma(b+q-a)(q-1)!(k-q)!}(b+q+z)_{k-q}(b+1-a+z)_{q-1}
$$
is a polynomial of degree $k-1$, and
\begin{equation}\label{eq:Bl111}
\sigma_{ql}=\frac{(-1)^{l}}{l!}{}F\!\!\left(\begin{matrix}-l,
b+q-1,\f+\m\\ b+q,\f\end{matrix}\right)
-\frac{(-1)^{l}(\f-b-q+1)_\m}{(b+q)_l(\f)_\m}.
\end{equation}
\end{theorem}
\textbf{Proof.} The proof goes along the same lines as the proof of Theorem~\ref{th:Degenerate1},
but with transformation \eqref{eq:Kar112} replaced by the transformation \cite[Theorem~4]{KarpPril1}
\begin{multline}\label{eq:Kar12-1}
F\left.\!\!\left(\begin{matrix}a, b, \f+\m\\b+k,\f\end{matrix}\right\vert x\right)
\!=\!\frac{\Gamma(b-a+1)}{\Gamma(b)(\f)_\m}(1-x)^{1-a}
F\left.\!\left(\!\begin{array}{l:}1,b+1-a \\b+k \end{array}\:T_{k-1}\:\right\vert x\right)\\
+\sum\limits_{q=1}^k\frac{(-1)^{q-1}(b)_{k}}{(b+q-1)(q-1)!(k-q)!}
\sum\limits_{l=0}^{m-1}\sigma_{ql}(a)_lx^l(1-x)^{-a-l}.~~~~\square
\end{multline}

\noindent\textbf{Remark.} Assuming $\lambda=-a,$ $vk=l,$ $uk=l$ in Lemma~\ref{lm:IPDsummation},
we obtain the summation formula
\begin{multline}
F\left.\!\left(\begin{matrix}a+l,d+l,\hh+\pp+l\\e+l,\hh+l\end{matrix}\:\right.\right)
=\frac{(-1)^{l}(\hh-d)_{\pp}\Gamma(e-a-d)\Gamma(1+d-e+a)\Gamma(e+l)}{(\hh+l)_{\pp}\Gamma(e-a)\Gamma(e-d)
\Gamma(l+d-e+a+p+1)}
\\
\times\sum_{j=0}^{p}\frac{(d-e+1)_{j}}{j!}
F\!\left(\begin{matrix}-j,1-\hh+d\\1-\hh+d-\pp\end{matrix}\right)(l+d)_{j}(l+d-e+a+j+1)_{p-j}.
\end{multline}
Thus, the second terms in equalities \eqref{eq:Kartransfor1} and  \eqref{eq:Kar2transform} are finite sums.

\begin{theorem}\label{th:Degenerate3}
Suppose $(c-g-m+1)_{m-1}\ne0$, $(c-a-m+1)_{m-1}\ne0$, $(1+a+g-c)_{m-1}\ne0$ and $\lambda=c-a-g-m+1$.
Assuming convergence of the hypergeometric functions involved we have the transformation
\begin{multline}\label{eq:Deqvect}
F\left.\!\!\left(\!\begin{matrix}a,g,d,b,\f+\m,\hh+\pp
\\
c,e,b+1,\f,\hh\end{matrix}\right.\right)
=\frac{(\f-b)_\m}{(\f)_\m}F\left.\!\!\left(\!\begin{matrix}a,g,d,b,\hh+\pp\\c,e,b+1,\hh\end{matrix}\right.\right)
\\
+\frac{\Gamma(e)\Gamma(e+\lambda-d)[(\f)_{\m}-(\f-b)_\m]}{\Gamma(e+\lambda)(d-e-\lambda+1)_p(\f)_\m(\hh)_\pp}\:
F\left.\!\!\left(\!\!\begin{array}{c:} d,g+\lambda,a+\lambda,\hat\llambda+1\\c,e+\lambda,\hat\llambda   \end{array}\:Y_p(1,0)\:\right.\!\!\right),
\end{multline}
where the polynomial $Y_p(1,0)$ is given in \eqref{eq:Ypolynomial} with $\lambda=c-a-g-m+1$ and
$\hat\llambda$ is the vector of zeros of the polynomial of degree $m-1$ defined by
\begin{equation}\label{eq:Lmhat}
\hat{L}_{m-1}(t)=\sum\limits_{l=0}^{m-1}\frac{(-1)^l\beta_l(a)_l(g)_l(t)_l}{(c-a-m+1)_l(c-g-m+1)_l}
{}_{3}F_{2}\!\left(\begin{matrix}-m+1+l,t+l,c-a-g-m+1\\c-a-m+1+l,c-g-m+1+l\end{matrix}\right)
\end{equation}
with $\beta_l$ from \eqref{eq:Bl1}.
\end{theorem}
\textbf{Proof.} The proof repeats that of Theorem~\ref{th:Degenerate1} with transformation \eqref{eq:Kar112} replaced by the transformation \cite[Theorem 5]{KarpPril1}
\begin{multline}\label{eq:Deqvect1}
F\left.\!\!\left(\!\begin{matrix}a,g,b,\f+\m\\c,
b+1,\f\end{matrix}\right\vert x\right)
=\frac{(\f-b)_\m}{(\f)_\m}{}_{3}F_{2}\left.\!\!\left(\!\begin{matrix}a,g,b\\c,b+1\end{matrix}\right\vert
x\right)
\\
+\frac{(\f)_{\m}-(\f-b)_\m}{(\f)_\m}(1-x)^{c-a-g-m+1}
F\left.\!\!\left(\!\begin{matrix}c-a-m+1,c-g-m+1,\hat\llambda+1\\c,\hat\llambda\end{matrix}\right\vert
x\right).~~~~~\square
\end{multline}

\section{Alternative approach: interchange of the order of summations}

In some situations we can exchange the order of summations on the right hand side of \eqref{eq:master1} to get the  hypergeometric function with several parameters shifted by unity as defined in \eqref{eq:pFqdefined}.
To illustrate this is idea we apply it to Euler's transformation
\begin{equation}\label{eq:Euler-2}
{}_{2}F_{1}\left.\!\!\left(\begin{matrix}a, b
\\c\end{matrix}\right\vert x\right)
=(1-x)^{c-a-b}{}_{2}F_{1}\left.\!\!\left(\begin{matrix}c-a,
c-b\\c\end{matrix}\right\vert x\right).
\end{equation}
This leads to
\begin{theorem}\label{th:Euler2-interchange}
Suppose $q\in\N$, $\a,\b\in\C^{p}$.
Assuming convergence of the hypergeometric functions involved we have
\begin{equation}\label{eq:transformation1}
F\left.\!\!\left(\begin{matrix}\gamma+q,\beta,\a\\\gamma,\b\end{matrix}\right.\right)
=F\left.\!\left(\!\begin{array}{l:}\beta+q,\a\\\b+q\end{array}\:P_M\:\right.\!\!\right),
\end{equation}
where the polynomial $P_{M}(x)$ of degree $M=qp$ is defined by
\begin{equation}\label{eq:Pqp}
P_{M}(x)=(\b+x)_q\cdot{}_{p+2}F_{p+1}\left.\!\!\left(\begin{matrix}-q,\gamma-\beta,\a+x
\\
\gamma,\b+x\end{matrix}\right.\right).
\end{equation}
\end{theorem}
\textbf{Remark.} Using relation \eqref{eq:FP-roots} formula \eqref{eq:transformation1} can be written in a more traditional notation as
$$
{}_{p+2}F_{p+1}\left.\!\!\left(\begin{matrix}\gamma+q,\beta,\a\\\gamma,\b\end{matrix}\right.\right)
={}_{p+2}F_{p+1}\left.\!\!\left(\begin{matrix}-q,\gamma-\beta,\a\\\gamma,\b\end{matrix}\right.\right)
{}_{M+p+1}F_{M+p}\left.\!\!\left(\begin{matrix}\beta+q,\a,1-\boldsymbol{\lambda}\\\b+q,-\boldsymbol{\lambda}\end{matrix}\right.\right),
$$
where $\boldsymbol{\lambda}=(\lambda_1,\ldots,\lambda_{qp})$ is the vector of roots of the polynomial \eqref{eq:Pqp}.

\textbf{Proof.} Application of Lemma~\ref{lm:master1} to Euler's transformation \eqref{eq:Euler-2} yields:
$$
{}_{p+2}F_{p+1}\left.\!\!\left(\begin{matrix}\alpha,\beta,\a\\\gamma,\b\end{matrix}\right.\right)
=\sum\limits_{k=0}^{\infty}\frac{(\gamma-\alpha)_k(\gamma-\beta)_k(\a)_k}{(\gamma)_k(\b)_kk!}
{}_{p+1}F_{p}\left.\!\!\left(\begin{matrix}\alpha+\beta-\gamma,\a+k\\\b+k\end{matrix}\right.\right).
$$
Now, assume that $\alpha=\gamma+q$, $q\in\N$.  Then exchanging the order of summations we get:
\begin{multline*}
{}_{p+2}F_{p+1}\left.\!\!\left(\begin{matrix}\gamma+q,\beta,\a\\\gamma,\b\end{matrix}\right.\right)
=\sum\limits_{k=0}^{q}\frac{(-q)_k(\gamma-\beta)_k(\a)_k}{(\gamma)_k(\b)_kk!}
{}_{p+1}F_{p}\left.\!\!\left(\begin{matrix}\beta+q,\a+k\\\b+k\end{matrix}\right.\right)
\\
=\sum\limits_{n=0}^{\infty}\frac{(\beta+q)_n}{n!}\sum\limits_{k=0}^{q}\frac{(-q)_k(\gamma-\beta)_k(\a)_k(\a+k)_n}
{(\gamma)_k(\b)_k(\b+k)_nk!}
\\
=\sum\limits_{n=0}^{\infty}\frac{(\beta+q)_n(\a)_n}{(\b)_nn!}\sum\limits_{k=0}^{q}\frac{(-q)_k(\gamma-\beta)_k(\a+n)_k}
{(\gamma)_k(\b+n)_kk!}
=\sum\limits_{n=0}^{\infty}\frac{(\beta+q)_n(\a)_n}{(\b)_n(\b+n)_qn!}P_{M}(n),
\end{multline*}
where $M=qp$ and $P_{M}(x)$ is defined in \eqref{eq:Pqp}.  This proves \eqref{eq:transformation1}.  To justify the expression from the remark denote the zeros of this polynomial by $\lambda_1$, $\lambda_2$, $\ldots$, $\lambda_{M}$ and note that the constant term of this polynomial, $P_{M}(0)$, is easily computed, so that in view of \eqref{eq:Pfactorized} we have 
$$
P_{M}(n)=(\b)_{q}\cdot{}_{p+2}F_{p+1}\left.\!\!\left(\begin{matrix}-q,\gamma-\beta,\a\\\gamma,\b\end{matrix}\right.\right)
\frac{(1-\lambda_1)_n\cdots(1-\lambda_{M})_n}{(-\lambda_1)_n\cdots(-\lambda_{M})_n}.
$$
It remains to apply $(\b)_n(\b+n)_q=(\b)_q(\b+q)_n$.$\hfill\square$

In particular, for $p=q=1$:
$$
{}_{3}F_{2}\left.\!\!\left(\begin{matrix}\gamma+1,\beta,a\\\gamma,b\end{matrix}\right\vert
1\right)
=\left(1-\frac{(\gamma-\beta)a}{\gamma{b}}\right){}_{3}F_{2}\left.\!\!\left(\begin{matrix}\beta+1,a,1-\lambda\\b+1,-\lambda\end{matrix}\right\vert
1\right),
$$
where $\lambda=\beta^{-1}((\gamma-\beta)a-\gamma{b})$.  For $m=1$,
$p=2$:
$$
{}_{4}F_{3}\left.\!\!\left(\begin{matrix}\gamma+1,\beta,a_1,a_2\\\gamma,b_1,b_2\end{matrix}\right\vert
1\right)
=\left(1-\frac{(\gamma-\beta)a_1a_2}{\gamma{b_1}b_2}\right)
{}_{5}F_{4}\left.\!\!\left(\begin{matrix}\beta+1,a_1,a_2,1-\lambda_1,1-\lambda_2\\b_1+1,b_2+1,-\lambda_1,-\lambda_2\end{matrix}\right\vert
1\right),
$$
where $\lambda_1$, $\lambda_2$ are the roots of
$$
P_2(x)=(a_1+x)(a_2+x)(\gamma-\beta)-\gamma(b_1+x)(b_2+x)=0.
$$

If we start with the first Euler-Pfaff transformation
\begin{equation}\label{eq:Euler-1}
{}_{2}F_{1}\left.\!\!\left(\begin{matrix}a, b
\\c\end{matrix}\right\vert x\right)
=(1-x)^{-a}{}_{2}F_{1}\left.\!\!\left(\begin{matrix}a, c-b
\\c\end{matrix}\right\vert\frac{x}{x-1}\right),
\end{equation}
we arrive at the following theorem.
\begin{theorem}\label{th:Euler1-interchange}
Suppose $-a_1=m\in\N$, $\N\ni{q}\le{m}$, $\a,\b\in\C^{p}$. Then
\begin{equation}\label{eq:transformation2}
F\left.\!\!\left(\!\begin{matrix}\gamma+q,\alpha,\a\\\gamma,\b\end{matrix}\right.\right)
=F\left.\!\left(\!\begin{array}{l:}\alpha,\a\\\b+q\end{array}\:R_{M}\:\right.\!\!\!\right),
\end{equation}
 where the polynomial $R_{M}$ of degree $M=q(p+1)$ is defined by
\begin{equation}\label{eq:polyRM}
R_{M}(x)=(\b+x)_{q}\sum\limits_{k=0}^{q}\frac{(-q)_k(\alpha+x)_{k}(\a+x)_{k}}{(\b+x)_{k}(-1)^k(\gamma)_kk!}
=(\b+x)_q\cdot{}_{p+2}F_{p+1}\!\!\left(\left.\begin{matrix}-q,\alpha+x,\a+x\\\gamma,\b+x\end{matrix}\right|-1\right).
\end{equation}
\end{theorem}
\textbf{Remark.}  Using relation \eqref{eq:FP-roots} formula \eqref{eq:transformation2} can be written in a more traditional notation as
\begin{equation}\label{eq:transformation2-2}
{}_{p+2}F_{p+1}\left.\!\!\left(\begin{matrix}\gamma+q,\alpha,\a\\\gamma,\b\end{matrix}\right.\right)
={}_{p+2}F_{p+1}\!\!\left(\left.\begin{matrix}-q,\alpha,\a\\\gamma,\b\end{matrix}\right|-1\right)
{}_{M+p+1}F_{M+p}\left.\!\!\left(\begin{matrix}\alpha,\a,1-\boldsymbol{\eta}\\\b+q,-\boldsymbol{\eta}\end{matrix}\right.\right),
\end{equation}
where $\boldsymbol{\eta}=(\eta_1,\ldots,\eta_{M})$ is the vector of roots of the polynomial \eqref{eq:polyRM}.

\textbf{Proof.}  Set $\beta=\gamma+q$. Application of Lemma~\ref{lm:master1} to Euler's transformation \eqref{eq:Euler-1} yields:
 \begin{multline*}
{}_{p+2}F_{p+1}\left.\!\!\left(\begin{matrix}\alpha,\beta,\a\\\gamma,\b\end{matrix}\right.\right)
=\sum\limits_{k=0}^{q}\frac{(\alpha)_k(\gamma-\beta)_k(\a)_k}{(-1)^k(\gamma)_k(\b)_kk!}
{}_{p+1}F_{p}\left.\!\!\left(\begin{matrix}\alpha+k,\a+k\\\b+k\end{matrix}\right.\right)
\\
=\sum\limits_{k=0}^{q}\frac{(\alpha)_k(-q)_k(\a)_k}{(-1)^k(\gamma)_k(\b)_kk!}
\sum\limits_{n=0}^{m}\frac{(\alpha+k)_{n}(\a+k)_{n}}{(\b+k)_{n}n!}
=\sum\limits_{n=0}^{m}\frac{(\alpha)_n(\a)_n}{(\b+n)_{q}(\b)_nn!}
R_{M}(n),
\end{multline*}
where $M=q(p+1)$ and $R_{M}(x)$ is defined in \eqref{eq:polyRM} which proves \eqref{eq:transformation2}.
Denoting the zeros of this polynomials by $\eta_1$, $\eta_2$, $\ldots$, $\eta_{M}$ and noting that the constant term of this polynomial, $R_{M}(0)$, is easily computed, we get
$$
R_{M}(n)=(\b)_q\cdot{}_{p+2}F_{p+1}\!\!\left(\left.\begin{matrix}-q,\alpha,\a\\\gamma,\b\end{matrix}\right|-1\right)
\frac{(1-\eta_1)_n\cdots(1-\eta_{M})_n}{(-\eta_1)_n\cdots(-\eta_{M})_n}.
$$
Substituting this expression back, we obtain \eqref{eq:transformation2-2}. $\hfill\square$

In particular, for $p=q=1$:
$$
{}_{3}F_{2}\left.\!\!\left(\begin{matrix}\gamma+1,\alpha,-m\\\gamma,b\end{matrix}\right\vert
1\right)
=\left(1-\frac{\alpha{m}}{\gamma{b}}\right){}_{4}F_{3}\left.\!\!\left(\begin{matrix}\alpha,-m,1-\eta_1,1-\eta_2\\b+1,-\eta_1,-\eta_2\end{matrix}\right\vert
1\right),
$$
where $\eta_1$, $\eta_2$ are the roots of
$$
\gamma(b+x)+(\alpha+x)(x-m)=0.
$$

This approach works for a number of other transformations listed in \cite[section~2.3]{KPGFmethod}. For another example we take Maier's recent generalization of Whipple's quadratic transformation for ${}_3F_2$ \cite[Theorem~3.1]{Maier2019}:
\begin{equation}\label{eq:MaierTh3.1}
F\!\left(\left.\!\begin{matrix}\alpha,\beta,\delta,1-\rrho\\1+\alpha-\beta,1+\alpha-\delta,-\rrho\end{matrix}\:\right|x\right)
=(1-x)^{-\alpha}{}_{3}F_{2}\!\left(\left.\!\begin{matrix}\alpha/2,\alpha/2+1/2,\alpha-\beta-\delta-r+1
\\
1+\alpha-\beta,1+\alpha-\delta\end{matrix}\right|\frac{-4x}{(1-x)^2}\right),
\end{equation}
where $\rrho$ is the vector of roots of the $2r$ degree polynomial
\begin{equation}\label{eq:MaierQ2kpoly}
P_{2r}(t;\alpha,\beta,\delta)={}_{3}F_{2}\!\left(\!\begin{matrix}-r,-t,t+\alpha\\\beta,\delta\end{matrix}\right).
\end{equation}

\begin{theorem}\label{th:Maier-interchange}
Suppose $r\in\N$ and $q\in\N$ satisfies $2q<\sum_{j=1}^{p}(b_j-a_j)-\alpha$. Then
\begin{equation}\label{eq:Maier-interchange}
F\left.\!\left(\!\!\begin{array}{l:}\beta+\delta+r-q-1,\beta,\delta,\a\\\delta+r-q,\beta+r-q,\b\end{array}\:P_{2r}\:\right.\!\!\right)
=\frac{1}{(\b)_{q}}
F\left.\!\!\left(\!\!\begin{array}{l:}\beta+\delta+r-q-1,\a\\\b+q\end{array}\:Q_{M}\:\right.\!\!\right),
\end{equation}
where $Q_M(x)$ is a polynomial of degree $M=(p+2)q$ given by
$$
Q_M(x)=(\b+x)_q\cdot
F\!\left(\left.\!\begin{matrix}-q,\Delta(\beta+\delta+r-q-1+x,2),\a+x\\\delta+r-q,\beta+r-q,\b+x\end{matrix}\:\right|-4\right).
$$
Formula \eqref{eq:Maier-interchange} remains true for $r=0$ if we omit $P_{2r}$ on the left hand side.
\end{theorem}
\textbf{Proof.}  Set $\alpha=\beta+\delta+r-q-1$.  Application of Lemma~\ref{lm:master1} to formula \eqref{eq:MaierTh3.1} then yields:
\begin{multline*}
F\left.\!\left(\!\!\begin{array}{l:}\alpha,\beta,\delta,\a\\\delta+r-q,\beta+r-q,\b\end{array}\:P_{2r}\:\right.\!\!\right)
=
\sum\limits_{k=0}^{q}\frac{\Delta(\alpha,2)_{k}(-q)_k(\a)_k(-4)^k}{(\delta+r-q)_{k}(\beta+r-q)_{k}(\b)_kk!}
F\!\left(\left.\!\!\begin{matrix}\alpha+2k,\a+k\\\b+k\end{matrix}\:\right.\!\!\right)
\\
=\sum\limits_{n=0}^{\infty}\frac{(\alpha)_n(\a)_{n}}{(\b)_n(\b+n)_{q}n!}
\sum\limits_{k=0}^{q}\frac{(-q)_k(\alpha+n)_{2k}(\a+n)_k(-1)^k(\b+n)_{q}}{(\delta+r-q)_{k}(\beta+r-q)_{k}(\b+n)_kk!}
=\sum\limits_{n=0}^{\infty}\frac{(\alpha)_n(\a)_{n}}{(\b)_n(\b+n)_{q}n!}Q_M(n),
\end{multline*}
where we applied the relations
$$
(\a)_{k}(\a+k)_{n}=(\a)_{n}(\a+n)_{k},~~(\alpha+n)_{2k}=4^{k}\Delta(\alpha+n,2)_{k}.~~~~~~\square
$$

\section{Summation formulas}
In this section we specialize some transformations from Section~2 and from our paper \cite{KPGFmethod} to get summation formulas which appeared interesting and new to us.  Note that the formulas presented in Theorems~\ref{th:GScubic-IPD-sum}, \ref{th:GScubicIPDsum}, \ref{th:GScubicSaal} sum hypergeometric functions with non-linearly constrained parameters.  This type of formulas is rarely met in the hypergeometric literature. Two examples were found by us recently  in \cite[(45)]{KPResults2019}, \cite[p.15]{KPGFmethod}.
\begin{theorem}\label{th:GScubic-IPD-sum}
For $\Re(\beta)>0$ the following summation formula holds true\emph{:}
\begin{equation}
{}_{5}F_{4}\left.\!\!\left(\!\begin{matrix}
\beta,d,1-d,\zeta+1,\psi+1
\\
\gamma,3+2\beta-\gamma,\zeta,\psi\end{matrix}\right.\!\!\right)
\!=\!\frac{A\Gamma(\gamma)\Gamma\left((d+\gamma)/2\right)\Gamma\left(1+\beta-(d+\gamma)/2\right)
\Gamma\left(3+2\beta-\gamma\right)}{\Gamma(d+\gamma-1)\Gamma\left((\gamma-d)/2\right)
\Gamma\left(1+\beta+(d-\gamma)/2\right)\Gamma(2+2\beta-d-\gamma)},
\end{equation}
where $\zeta=(2+\beta-d-\gamma)(\gamma-d-1)/\beta+\gamma-1$,
$$
\psi=\frac{\beta(1-d)(\beta+d)}{(\beta+d)(2+\beta-d-\gamma)+(\gamma-\beta-1)(d+\gamma-2)}
$$
and
$$
A=\frac{\beta}{(\beta(\gamma-1)+(2+\beta-d-\gamma)(\gamma-d-1))(\beta+d)}.
$$
\end{theorem}
\textbf{Proof.} In \cite[(59)]{KPGFmethod} we proved that
\begin{multline}\label{eq:5F4-2unitshifts}
{}_{5}F_{4}\left.\!\!\left(\!\begin{matrix}
a,b,d,h+1,f+1\\c,e,h,f\end{matrix}\right.\right)
=
\frac{((e-d-1)h+(c-a-b-1)(h-d))\Gamma(e)\Gamma(s_*)}{h\Gamma(s_*+d+1)\Gamma(e-d)}
\\
\times{}_{5}F_{4}\left.\!\!\left(\!\begin{matrix}
c-a-1,c-b-1,d,\hat{\xi}+1,\hat{\zeta}+1\\c,e+c-a-b-1,\hat{\xi},\hat{\zeta}\end{matrix}\right.\right),
\end{multline}
where $s_*=e+c-a-b-d-2$, $c-a-1\ne0$, $c-b-1\ne0$,
$$
\hat{\xi}=h+\frac{(c-a-b-1)(h-d)}{e-d-1}~~\text{and}~~\hat{\zeta}=\frac{(c-a-1)(c-b-1)f}{(c-a-b-1)f+ab}.
$$
Setting $e=f+1$, $c=h+1$, we will have
\begin{equation}
{}_{3}F_{2}\left.\!\!\left(\!\begin{matrix}
a,b,d\\h,f\end{matrix}\right.\right)
=
L\cdot{}_{5}F_{4}\left.\!\!\left(\!\begin{matrix}
h-a,h-b,d,\psi+1,\zeta+1\\h+1,f+h-a-b+1,\psi,\zeta\end{matrix}\right.\right),
\end{equation}
where
$$
\zeta=h+\frac{(h-a-b)(h-d)}{(f-d)},~~~\psi=\frac{(h-a)(h-b)f}{(h-a-b)f+ab},
$$
$$
L=\frac{((f-d)h+(h-a-b)(h-d))\Gamma(f+ 1))\Gamma(h+f-a-b-d)}{(f-d)h\Gamma(h+f-a-b+1)\Gamma(f-d)}.
$$
Assume that $f=d+\beta$, $h=\gamma-1$, $a=\gamma+d-2$, $b=\gamma-\beta-1$. Then we can apply Dixon's theorem (see, for example, \cite[(2.2.11)]{AAR}) to sum
$$
{ }_{3}F_{2}\left.\!\!\left(\!\begin{matrix}
\gamma+d-2,\gamma-\beta-1,d\\d+\beta,\gamma-1\end{matrix}\right.\right),
$$
which after some algebra yields the result. $\hfill\square$

Next theorem is a summation formula for general very well-poised  non-terminating ${}_{7}F_{6}$ containing a parameter pair $\left[\!\!\begin{array}{l}F+1\\F\end{array}\!\!\right]$.
\begin{theorem}\label{th:7F6unitshift}
Suppose $\Re(A-C-D-E)>0$. Then
\begin{multline}\label{eq:7F6nonlinsum}
{}_{7}F_{6}\left.\!\!\left(\!\begin{matrix} A,1+A/2,C,D,E,F+1,1+A-F
\\
A/2,1+A-C,1+A-D,1+A-E,F,A-F\end{matrix}\right.\right)
\\
=\frac{\Gamma(1+A-C)\Gamma(1+A-D)\Gamma(1+A-E)\Gamma(1+A-C-D-E)}{\Gamma(1+A)\Gamma(1+A-C-E)\Gamma(1+A-D-E)\Gamma(1+A-C-D)}
\!\left(\!1+\frac{CDE}{F(F-A)(C+D+E-A)}\!\right).
\end{multline}
\end{theorem}
\textbf{Proof.}  Changing capital to lowercase letters in \cite[(79)]{KPGFmethod} and setting $(a-1)/2-b=-1$ we will have:
\begin{multline*}
{ }_{7}F_{6}\left.\!\!\left(\!\begin{matrix} 2b-1,b+1/2,c,d,e,b+1/2-\sigma,b+1/2+\sigma
\\
b-1/2,2b-c,2b-d,2b-e,b-1/2-\sigma,b-1/2+\sigma\end{matrix}\right.\right)
\\
=\frac{\Gamma(2b-c)\Gamma(2b-d)\Gamma(2b-e)\Gamma(2b-c-d-e)}{\Gamma(2b-c-e)\Gamma(2b-d-e)\Gamma(2b-c-d)\Gamma(2b)}
\!\left(1-\frac{cde(f+1)}{b^2f(c+d+e-2b+1)}\right),
\end{multline*}
where
$$
\sigma^2=\frac{(2b-1)^2-(4b-1)f}{4(f+1)}.
$$
The last expression can be written as
$$
f=\frac{(b-1/2)^2-\sigma^2}{(b-1/4)+\sigma^2}.
$$
Using this expression we can get rid of $f$ on the  right hand side of the last formula.
Next, setting $A=b-1/2$, $F=b-1/2-\sigma$, $C=c$, $D=d$, $E=e$,  after much rearrangements and simplifications we arrive at  \eqref{eq:7F6nonlinsum}.$\hfill\square$

\begin{theorem}\label{th:GScubicIPDsum}
For $-d\in\N$ the following summation formula holds:
\begin{multline}\label{eq:GScubicIPDsum}
{}_7F_{6}\left.\!\!\left(\begin{matrix}\Delta(\alpha,3),\Delta(1+\alpha-e,2),d,\xi+1\\
\alpha/2,(\alpha+3)/2,\Delta(2+\alpha+d-e,3),\xi\end{matrix}\right.\right)
\\
=\frac{h\Gamma(e-\alpha)\Gamma(e-d)(1+\alpha+d-e)}{\Gamma(e)\Gamma(e-\alpha-d)((h-d)\alpha-(e-d-1)h)}\left(1-\frac{2d(h+1)}{(\alpha+3)eh}\right),
\end{multline}
where
$$
\xi=\frac{(h-d)\alpha-(e-d-1)h}{3h-2d-e+1}.
$$
\end{theorem}
\textbf{Proof.} Set $2\beta=\alpha+1$ in \eqref{eq:GScubic-IPD}, use \eqref{eq:Y1} for $Y_1(1,3)$ and simplify.$\hfill\square$

\begin{theorem}\label{th:GScubicSaal}
For $-d\in\N$ the following summation formula holds:
\begin{multline}\label{eq:GScubicSaal}
{}_8F_{7}\left.\!\!\left(\begin{matrix}\Delta(\alpha,3),\beta,\Delta(1+\alpha-e,2),d,\xi+1\\
\Delta(\alpha,2),\beta+1,\Delta(2+\alpha+d-e,3),\xi\end{matrix}\right.\right)
\\
=\frac{\Gamma(e-\alpha)\Gamma(e-d-1)\beta(1+\alpha+d-e)[h(3\beta-\alpha)(\beta-d-1)+d(2\beta-\alpha+h)]}{\Gamma(e-1)\Gamma(e-\alpha-d)(\beta-d-1)(\beta-d)(3\beta-\alpha)((h-d)\alpha-(e-d-1)h)},
\end{multline}
where $3\beta+e-\alpha-d=2$ and
$$
\xi=\frac{(h-d)\alpha-(e-d-1)h}{3h-2d-e+1}.
$$
\end{theorem}
\textbf{Proof.}  Rakha and Rathie \cite[(2.5)]{RakhaRathie} (see also \cite[(3.1)]{KRP2013})
extended the Pfaff-Saalsch\"{u}tz summation theorem by adding a parameter pair with unit shift. Their extension can be written in the form:
$$
{}_{4}F_{3}\!\left(\begin{matrix}-j,a,b,f+1\\c,d,f\end{matrix}\right)
=\frac{(c-a-1)_{j}(d-a)_{j}(\gamma+1)_{j}}{(c)_{j}(d)_{j}(\gamma)_{j}},~~~~c+d=-j+a+b+2,
$$
where
$$
\gamma=\frac{(c-a-1)(c-b-1)f}{ab+(c-a-b-1)f}.
$$
Imposing the condition $3\beta+e-\alpha-d=2$ we get Saalsch\"{u}tzian ${}_{4}F_{3}$ with one unit shift on the left hand side of \eqref{eq:GScubic-IPD}. Using the above formula to sum  ${}_{4}F_{3}$  and applying \eqref{eq:Y1} for $Y_1(1,3)$ after some simplifications we arrive at  \eqref{eq:GScubicSaal}. $\hfill\square$

\begin{theorem}\label{th:5F4sum}
For each $r,n\in\N$ the following summation formula is true:
\begin{multline}\label{eq:5F4sum}
{}_{5}F_{4}\left.\!\!\left(\begin{matrix}\Delta(\alpha-1,3),-1/2-n,-n\\
\alpha/2+r,\alpha/2+3/2,-1/3-n,-2/3-n\end{matrix}\right.\right)
=
\frac{3\Gamma(2n+3)\Gamma(3n+\alpha+3)}{(\alpha+2)\Gamma(3n+4)\Gamma(2n+\alpha+2)}
\\
\times\left\{1+\frac{2n\alpha(\alpha-1)(1-r)(\beta-r-1/2)_{r}(-\beta-r-1/2)_{r}}{(\alpha+3)(\alpha+5)(\alpha+2r)(\alpha+2n+2)(\beta+1/2)_{r}(-\beta+1/2)_{r}}\left[1+\frac{2r(\alpha+1)}{(\beta-r-1/2)(-\beta-r-1/2)}\right]\right\}.
\end{multline}
\end{theorem}
\textbf{Proof}. Set $\beta=3/2-r$ in \eqref{eq:Maier3.2-GS}.$\hfill\square$

For $r=0$ the above theorem takes a particularly simple form:
\begin{multline*}
{}_{5}F_{4}\left.\!\!\left(\begin{matrix}\Delta(\alpha-1,3),-1/2-n,-n\\
\alpha/2,\alpha/2+3/2,-1/3-n,-2/3-n\end{matrix}\right.\right)
\\
=
\frac{3\Gamma(2n+3)\Gamma(3n+\alpha+3)}{(\alpha+2)\Gamma(3n+4)\Gamma(2n+\alpha+2)}\left(1+\frac{2n(\alpha-1)}{(\alpha+3)(\alpha+5)(\alpha+2n+2)}\right).
\end{multline*}

\textbf{Acknowledgements.}  The second author has been supported by the Russian Foundation for Basic Research under project 19-010-00206.

\end{document}